\documentclass[a4paper,11pt]{amsart}

\usepackage{graphicx}
\usepackage{mathptmx}
\usepackage{amsmath}
\usepackage{amssymb}
\usepackage{enumitem}
\usepackage{xcolor}

\newmuskip\pFqmuskip

\newcommand*\pFq[6][8]{%
  \begingroup 
  \pFqmuskip=#1mu\relax
  \mathcode`=\string"8000
  \begingroup\lccode`\~=`\,
  \lowercase{\endgroup\let~}\pFqcomma
  F^{#2}_{#3}{\left(\genfrac..{0pt}{}{#4}{#5}\bigg|#6\right)}%
  \endgroup
}
\newcommand{\pFqcomma}{\mskip\pFqmuskip}

\newtheorem{theorem}{Theorem}[section]

\newtheorem{proposition}[theorem]{Proposition}

\begin{document}

\title[A probabilistic proof of a recurrence relation]{A probabilistic proof of a recurrence relation for sums of values of degenerate falling factorials}

\author{Taekyun  Kim}
\address{Department of Mathematics, Kwangwoon University, Seoul 139-701, Republic of Korea}
\email{tkkim@kw.ac.kr}
\author{Dae San  Kim}
\address{Department of Mathematics, Sogang University, Seoul 121-742, Republic of Korea}
\email{dskim@sogang.ac.kr}


\begin{abstract}
In this paper, we consider sums of values of degenerate falling factorials and give a probabilistic proof of a recurrence relation for them. This may be viewed as a degenerate version of the recent probabilistic proofs on sums of powers of integers.
\end{abstract}

\maketitle

\section{Introduction} 
Jacob Bernoulli considered sums of powers of the first $n$ positive integers, $1^{k}+2^{k}+\cdots+n^{k}$, which have been a topic of research for centuries. We note that 
\begin{align*}
	&1+2+3+\cdots+n=\binom{n+1}{2},\\
	&1^{2}+2^{2}+3^{2}+\cdots+n^{2}=\frac{n(n+1)(2n+1)}{6},\\
	&1^{3}+2^{3}+3^{3}+\cdots+n^{3}=\big(1+2+3+\cdots+n\big)^{2},\quad (\mathrm{see}\ [6,7,9-11]).
\end{align*} \par
The Bernoulli polynomials are defined by  
\begin{equation}
\frac{t}{e^{t}-1}e^{xt}=\sum_{n=0}^{\infty}B_{n}(x)\frac{t^{n}}{n!},\quad (\mathrm{see}\ [1-20]). \label{1}
\end{equation}
When $x=0,\ B_{n}=B_{n}(0)$ are called the Bernoulli numbers. 
By \eqref{1}, we get 
\begin{equation}
B_{n}(x)=\sum_{k=0}^{n}\binom{n}{k}B_{k}x^{n-k},\quad (n\ge 0),\quad (\mathrm{see}\ [10,12]).\label{2}
\end{equation}
From \eqref{1}, we note that 
\begin{equation}
\sum_{k=0}^{n-1}(k+x)^{m}=\frac{1}{m+1}\sum_{l=0}^{m}\binom{m+1}{l}B_{l}(x)n^{m+1-l},\label{3}
\end{equation}
where $m,n\in\mathbb{N}$. \par 
Denoting the sum $1^{m}+2^{m}+\cdots+n^{m}$ by $S_{m}(n)$, by \eqref{3}, we get 
\begin{align}
S_{m}(n)&=\frac{1}{m+1}\sum_{k=0}^{m}\binom{m+1}{l}B_{l}(n+1)^{m+1-l}=\int_{0}^{n+1} B_{m}(u) du, \label{4}\\
S_{m}(n)&=\frac{(n+1)^{m+1}-1}{m+1}-\frac{1}{m+1}\sum_{r=0}^{m-1}\binom{m+1}{r}S_{r}(n), \label{4-1}\\
&=\frac{n^{m+1}}{m+1}+\sum_{r=0}^{m-1}\binom{m}{r}\frac{(-1)^{m-r+1}}{m-r+1}S_{r}(n),\label{4-2}
\end{align}
where $m$ is a positive integer (see [6,7,10]). \par 
For any $\lambda\in\mathbb{R}$, the degenerate exponentials are defined by 
\begin{equation}
e_{\lambda}^{x}(t)=\sum_{n=0}^{\infty}(x)_{n,\lambda}\frac{t^{n}}{n!},\quad (\mathrm{see}\ [11-15]),\label{5}
\end{equation}
where the degenerate falling factorials are given by 
\begin{equation*}
(x)_{0,\lambda}=1,\ (x)_{n,\lambda}=x(x-\lambda)(x-2\lambda)\cdots(x-(n-1)\lambda),\ (n\ge 1).
\end{equation*}
In particular, for $x=1$, we denote them by $e_{\lambda}(t)=e_{\lambda}^{1}(t)$. \par
The degenerate Bernoulli polynomials are defined by 
\begin{equation}
\frac{t}{e_{\lambda}(t)-1}e_{\lambda}^{x}(t)=\sum_{n=0}^{\infty}\beta_{n,\lambda}(x)\frac{t^{n}}{n!},\quad (\mathrm{see}\ [4,16,20]). \label{6}
\end{equation}
When $x=0$, $\beta_{n,\lambda}=\beta_{n,\lambda}(0)$ are called the degenerate Bernoulli numbers. \par 
The degenerate Stirling numbers of the second kind are defined by Kim-Kim as 
\begin{equation}
(x)_{n,\lambda}=\sum_{k=0}^{n}{n \brace k}_{\lambda}(x)_{k},\quad (n\ge 0),\quad (\mathrm{see}\ [8]), \label{7}
\end{equation}
where the falling factorials are given by
\begin{equation*}
(x)_{0}=1,\ (x)_{n}=x(x-1)(x-2)\cdots(x-n+1),\ (n\ge 1).
\end{equation*} \par
In this paper, we study sums of values of degenerate falling factorials which are given by
\begin{equation}
S_{k,\lambda}(n)=(1)_{k,\lambda}+(2)_{k,\lambda}+\cdots+(n)_{k,\lambda}=\sum_{j=1}^{n}(j)_{k,\lambda},\quad (k\in\mathbb{N}). \label{8}	
\end{equation}
In Section 1, we recall the necessary facts that are needed throughout this paper. After recalling two expressions of $S_{k,\lambda}(n)$, we derive two recurrence relations for them in Section 2. Let $X$ be a nonnegative integer-valued random variable, and let $k$ be a positive integer. Then, in Section 3, we show first that the $k$-th degenerate moment of $X$ is given by $E\Big[(X)_{k,\lambda}\Big]=\sum_{x=0}^{\infty}\big((x+1)_{k,\lambda}-(x)_{k,\lambda}\big)P\{X>x\}$. Then we apply this to the uniform random variable $X$ supported on $\{0,1,2,\dots , n\}$ to derive a recurrence relation for $S_{k,\lambda}(n)$. 

\section{Some formulas for sums of values of degenerate falling factorials} 
From \eqref{6}, we note that 
\begin{equation}
\sum_{m=0}^{\infty}\sum_{k=0}^{n-1}(k+x)_{m,\lambda}\frac{t^{m}}{m!}=\sum_{m=0}^{\infty}\frac{1}{m+1}\Big(\beta_{m+1,\lambda}(n+x)-\beta_{m+1,\lambda}(x)\Big)\frac{t^{m}}{m!}. \label{9}	
\end{equation}
Thus, by \eqref{9}, we get 
\begin{align}
\sum_{k=0}^{n-1}(k+x)_{m,\lambda}&=\frac{1}{m+1}\Big(\beta_{m+1,\lambda}(n+x)-\beta_{m+1,\lambda}(x)\Big)\label{10} \\
&=\frac{1}{m+1}\sum_{l=0}^{m}\binom{m+1}{l}\beta_{l,\lambda}(x)(n)_{m+1-l,\lambda}, \nonumber	
\end{align}
where $m$ is nonnegative integer. \\
From \eqref{10}, for $m, n \in\mathbb{N}$ and $x=0$, we have the next proposition. This is also obtained in [14, Lemma 7].
\begin{proposition}
For $m,n\in\mathbb{N}$, we have 
\begin{displaymath}
S_{m,\lambda}(n-1)=\frac{1}{m+1}\sum_{l=0}^{m}\binom{m+1}{l}(n)_{m+1-l,\lambda}\beta_{l,\lambda}.
\end{displaymath}
\end{proposition}

As it is done in [14, Theorem 8], by using \eqref{7} and \eqref{8}, we get 
the following theorem. 
\begin{theorem}
For $k\in\mathbb{N}$, we have 
\begin{displaymath}
S_{k,\lambda}(n)=\sum_{l=1}^{k}{k \brace l}_{\lambda}\binom{n+1}{l+1}l!.
\end{displaymath}
\end{theorem}
Now, we observe that 
\begin{equation}
(x+y)_{n,\lambda}=\sum_{k=0}^{n}\binom{n}{k}(x)_{k,\lambda}(y)_{n-k,\lambda},\quad (n\ge 0). \label{13}
\end{equation}
From \eqref{8}, we note that, for $k\in\mathbb{N}$,
\begin{align}
\sum_{r=0}^{k}\binom{k+1}{r}(1)_{k+1-r,\lambda}S_{r,\lambda}(n)&=\sum_{j=1}^{n}\sum_{r=0}^{k}\binom{k+1}{r}(1)_{k+1-r,\lambda}(j)_{r,\lambda}\label{14}\\
&=\sum_{j=1}^{n}\bigg(\sum_{r=0}^{k+1}\binom{k+1}{r}(1)_{k+1-r,\lambda}(j)_{r,\lambda}-(j)_{k+1,\lambda}\bigg)\nonumber\\
&=\sum_{j=1}^{n}\Big((j+1)_{k+1,\lambda}-(j)_{k+1,\lambda}\Big) \nonumber\\
&=(n+1)_{k+1,\lambda}-(1)_{k+1,\lambda}. \nonumber
\end{align}
By \eqref{14}, we get 
\begin{align}
(n+1)_{k+1,\lambda}-(1)_{k+1,\lambda}&=\sum_{r=0}^{k}\binom{k+1}{r}(1)_{k+1-r,\lambda}S_{r,\lambda}(n)\label{15}\\
&=\sum_{r=0}^{k-1}\binom{k+1}{r}(1)_{k+1-r,\lambda}S_{r,\lambda}(n)+(k+1)S_{k,\lambda}(n).\nonumber
\end{align}
Thus, by \eqref{15}, we obtain the following theorem. 
\begin{theorem}
For $k\in\mathbb{N}$, we have the recurrence relation
\begin{displaymath}
S_{k,\lambda}(n)=\frac{(n+1)_{k+1,\lambda}-(1)_{k+1,\lambda}}{k+1}-\frac{1}{k+1}\sum_{r=0}^{k-1}\binom{k+1}{r}(1)_{k+1-r,\lambda}S_{r,\lambda}(n).
\end{displaymath}
\end{theorem}
We note that 
\begin{equation}
(n)_{k+1,\lambda}=\sum_{j=1}^{n}\Big((j)_{k+1,\lambda}-(j-1)_{k+1,\lambda}\Big),\quad (k\in\mathbb{N}).\label{17}
\end{equation}
From \eqref{13}, we note that 
\begin{align}
&(j)_{k+1,\lambda}-(j-1)_{k+1,\lambda}=(j)_{k+1,\lambda}-\sum_{r=0}^{k+1}\binom{k+1}{r}(-1)_{k+1-r,\lambda}(j)_{r,\lambda}\label{18}\\
&=(j)_{k+1,\lambda}-\sum_{r=0}^{k+1}\binom{k+1}{r}(-1)^{k+1-r}\langle 1\rangle_{k+1-r,\lambda}(j)_{r,\lambda}\nonumber\\
&=(j)_{k+1,\lambda}-(j)_{k+1,\lambda}+(k+1)(j)_{k,\lambda}-\sum_{r=0}^{k-1}\binom{k+1}{r}(-1)^{k+1-r}\langle 1\rangle_{k+1-r,\lambda}(j)_{r,\lambda}\nonumber\\
&=(k+1)(j)_{k,\lambda}-\sum_{r=0}^{k-1}\binom{k+1}{r}(-1)^{k+1-r}\langle 1\rangle_{k+1-r,\lambda}(j)_{r,\lambda},\quad (k\in\mathbb{N}), \nonumber
\end{align}
where the degenerate rising factorials are given by
\begin{equation*}
\langle x\rangle_{0,\lambda}=1,\ \langle x\rangle_{k,\lambda}=x(x+\lambda)(x+2\lambda)\cdots(x+(k-1)\lambda),\ (k\ge 1).
\end{equation*}
Thus, by \eqref{17} and \eqref{18}, we get 
\begin{align}
(n)_{k+1,\lambda}&=\sum_{j=1}^{n}\Big((j)_{k+1,\lambda}-(j-1)_{k+1,\lambda}\Big) \label{19} \\
&=(k+1)\sum_{j=1}^{n}(j)_{k,\lambda}-\sum_{r=0}^{k-1}\binom{k+1}{r}(-1)^{k+1-r}\langle 1\rangle_{k+1-r,\lambda}\sum_{j=1}^{n}(j)_{r,\lambda}\nonumber\\
&=(k+1)S_{k,\lambda}(n)-\sum_{r=0}^{k-1}\binom{k+1}{r}(-1)^{k+1-r}\langle 1\rangle_{k+1-r,\lambda}S_{r,\lambda}(n).\nonumber	
\end{align}
From \eqref{19}, we obtain the following theorem. 
\begin{theorem}
For $k \in\mathbb{N}$, we have 
\begin{displaymath}
S_{k,\lambda}(n)=\frac{(n)_{k+1,\lambda}}{k+1}+\frac{1}{k+1}\sum_{r=0}^{k-1}\binom{k+1}{r}(-1)^{k+1-r}\langle 1\rangle_{k+1-r,\lambda}S_{r,\lambda}(n).
\end{displaymath}
\end{theorem}

\section{A probabilistic proof of a recurrence relation for $S_{k,\lambda}(n)$} 
Recently, probabilistic methods are used in deriving recurrence formulas for sums of powers of integers, (see [6,7]). In this section, we give a probabilistic proof of a recurrence relation for sums of values of degenerate falling factorials.\par 
 Let $X$ be a nonnegative integer-valued random variable, and let $k$ be any positive integer. Then we note that 
\begin{align}
&\sum_{x=0}^{\infty}\Big((x+1)_{k,\lambda}-(x)_{k,\lambda}\Big)P\{X>x\} \label{21} \\
&=(1)_{k,\lambda}P\{X>0\}+\Big((2)_{k,\lambda}-(1)_{k,\lambda}\Big)P\{X>1\}+\Big((3)_{k,\lambda}-(2)_{k,\lambda}\Big)P\{X>2\}\nonumber\\
&\quad+\Big((4)_{k,\lambda}-(3)_{k,\lambda}\Big)P\{X>3\}+\cdots\nonumber \\
&=(1)_{k,\lambda}P\{X=1\}+(1)_{k,\lambda}P\{X>1\}-(1)_{k,\lambda}P\{X>1\}+(2)_{k,\lambda}P\{X=2\}\nonumber \\
&\quad +(2)_{k,\lambda}P\{X>2\} -(2)_{k,\lambda}P\{X>2\}+(3)_{k,\lambda}P\{X=3\}-(3)_{k,\lambda}P\{X>3\}\nonumber\\
&=(1)_{k,\lambda}P\{X=1\}+(2)_{k,\lambda}P\{X=2\}+(3)_{k,\lambda}P\{X=3\}+\cdots\nonumber \\
&=\sum_{x=0}^{\infty}(x)_{k,\lambda}P\{X=x\}=E\Big[(X)_{k,\lambda}\Big]. \nonumber
\end{align}
Therefore, by \eqref{21}, we obtain the following theorem. 
\begin{theorem}
Let $X$ be a nonnegative integer-valued random variable. For $k\in\mathbb{N}$, the $k$-th degenerate moment of $X$ is given by 
\begin{displaymath}
E\Big[(X)_{k,\lambda}\Big]=\sum_{x=0}^{\infty}\Big((x+1)_{k,\lambda}-(x)_{k,\lambda}\Big)P\{X>x\}.
\end{displaymath}
\end{theorem}
Assume that $X$ has support in $\{0,1,2,\dots,n\}$. Then we have 
\begin{align}
\sum_{x=0}^{n}(x)_{k,\lambda}P\{X=x\}&=\sum_{x=0}^{n}\Big((x+1)_{k,\lambda}-(x)_{k,\lambda}\Big)P\{X>x\}\label{22}\\
&=\sum_{x=0}^{n}\Big((x+1)_{k,\lambda}-(x)_{k,\lambda}\Big)\sum_{y=x+1}^{n}P\{X=y\},\nonumber
\end{align}
where $k$ is a positive integer. \par 
Now, let $X$ be the uniform random variable supported on $\{0,1,2,\dots,n\}$, that is, $P\{X=x\}=\frac{1}{n+1}$, for $ x \in \{0,1,2,\dots,n\}$. Then we note that 
\begin{equation*}
\sum_{y=x+1}^{n}P\{X=y\}=\frac{n-x}{n+1}. 
\end{equation*}
From \eqref{22}, we note that 
\begin{align}
\sum_{x=0}^{n}(x)_{k,\lambda}\frac{1}{n+1}&=\sum_{x=0}^{n}(x)_{k,\lambda}P\{X=x\} \label{23} \\
&=\sum_{x=0}^{n}\Big((x+1)_{k,\lambda}-(x)_{k,\lambda}\Big)\sum_{y=x+1}^{n}P\{X=y\}\nonumber\\
&=\sum_{x=0}^{n}\Big((x+1)_{k,\lambda}-(x)_{k,\lambda}\Big)\frac{n-x}{n+1},\quad (k\in\mathbb{N}).\nonumber
\end{align}
By \eqref{23}, we get 
\begin{align}
&S_{k,\lambda}(n)=\sum_{x=0}^{n}(x)_{k,\lambda}=\sum_{x=0}^{n}\Big((x+1)_{k,\lambda}-(x)_{k,\lambda}\Big)(n-x)\label{24}\\
&=n\sum_{x=0}^{n}\Big((x+1)_{k,\lambda}-(x)_{k,\lambda}\Big)-\sum_{x=0}^{n}x\Big((x+1)_{k,\lambda}-(x)_{k,\lambda}\Big)\nonumber \\
&=n(n+1)_{k,\lambda}-\sum_{x=0}^{n}x\bigg(\sum_{r=0}^{k}\binom{k}{r}(x)_{r,\lambda}(1)_{k-r,\lambda}-(x)_{k,\lambda}\bigg) \nonumber \\
&=n(n+1)_{k,\lambda}-\sum_{x=0}^{n}\sum_{r=0}^{k-1}\binom{k}{r}(x-r\lambda+r\lambda)(x)_{r,\lambda}(1)_{k-r,\lambda}\nonumber\\
&=n(n+1)_{k,\lambda}-\sum_{x=0}^{n}\sum_{r=0}^{k-1}\binom{k}{r}(x)_{r+1,\lambda}(1)_{k-r,\lambda}\nonumber\\
&\qquad\qquad\qquad-\lambda\sum_{x=0}^{n}\sum_{r=0}^{k-1}\binom{k}{r}r(x)_{r,\lambda}(1)_{k-r,\lambda} \nonumber \\
&=n(n+1)_{k,\lambda}-\sum_{x=0}^{n}\sum_{r=0}^{k-2}\binom{k}{r}(x)_{r+1,\lambda}(1)_{k-r,\lambda}-k\sum_{x=0}^{n}(x)_{k,\lambda}\nonumber\\
&\qquad\qquad\qquad-\lambda\sum_{x=0}^{n}\sum_{r=0}^{k-1}r\binom{k}{r}(x)_{r,\lambda}(1)_{k-r,\lambda}\nonumber \\
&=n(n+1)_{k,\lambda}-\sum_{r=0}^{k-2}(1)_{k-r,\lambda}\binom{k}{r}S_{r+1,\lambda}(n)-kS_{k,\lambda}(n)\nonumber\\
&\qquad\qquad\qquad -\lambda\sum_{r=1}^{k-1}r\binom{k}{r}(1)_{k-r,\lambda}S_{r,\lambda}(n)\nonumber\\
&=n(n+1)_{k,\lambda}-\sum_{r=1}^{k-1}(1)_{k+1-r,\lambda}\binom{k}{r-1}S_{r,\lambda}(n)-kS_{k,\lambda}(n)\nonumber\\
&\qquad\qquad\qquad-\lambda\sum_{r=1}^{k-1}r\binom{k}{r}(1)_{k-r,\lambda}S_{r,\lambda}(n), \nonumber
\end{align}
where $k$ is a positive integer. \\
By \eqref{24}, we obtain the following theorem. 
\begin{theorem}
For $k\in\mathbb{N}$, we have 
\begin{align*}
S_{k,\lambda}(n)&=\frac{n(n+1)_{k,\lambda}}{k+1}-\frac{1}{k+1}\sum_{r=1}^{k-1}(1)_{k+1-r,\lambda}\binom{k}{r-1}
S_{r,\lambda}(n) \nonumber \\
&\qquad\qquad\qquad-\frac{\lambda}{k+1}\sum_{r=1}^{k-1}r\binom{k}{r}(1)_{k-r,\lambda}S_{r,\lambda}(n).
\nonumber
\end{align*}
\end{theorem}
\section{Conclusion} 
In this paper, we derived three recurrence relations for the sums of values of degenerate falling factorials $S_{k,\lambda}(n)=(1)_{k,\lambda}+(2)_{k,\lambda}+\cdots+(n)_{k,\lambda},\quad (k\in\mathbb{N})$. They are Theorems 2.3, 2.4 and 3.2. If we let $\lambda \rightarrow 0$, then Theorem 2.3 and Theorem 2.4 boil down to \eqref{4-1} and \eqref{4-2}, respectively. In addition, we obtain another recurrence relation for $S_{k}(n)$ by letting $\lambda \rightarrow 0$. Namely, we get 
\begin{equation*}
S_{k}(n)=\frac{n(n+1)^{k}}{k+1}-\frac{1}{k+1}\sum_{r=1}^{k-1}\binom{k}{r-1}
S_{r}(n).
\end{equation*}


\begin{thebibliography}{9}
\bibitem{1}
Aydin, M. S.; Acikgoz, M.; Araci, S. \emph{A new construction on the degenerate Hurwitz-zeta function associated with certain applications,} Proc. Jangjeon Math. Soc. \textbf{25} (2022), no. 2, 195-203.
\bibitem{2}
Barman, K.; Chakraborty, B.; Morthini, R. \emph{Two Classical formulas for the sums of powers of consecutive integers via complex analysis,} Complex Anal. Synerg. \textbf{2024} (2024), Article no. 10:5.
\bibitem{3}
Beardon, A. F. \emph{Sums of powers of integers,} Amer. Math. Monthly \textbf{103} (1996), no. 3, 201-213.
\bibitem{4}
Carlitz, L. \emph{Degenerate Stirling, Bernoulli and Eulerian numbers,} Utilitas Math. \textbf{15} (1979), 51-88. 
\bibitem{5}
Chen, L.; Dolgy, D. V.; Kim, T.; Kim, D. S. \emph{Probabilistic type 2 Bernoulli and Euler polynomials,} AIMS Math. \textbf{9} (2024), no. 6, 14312-14324.
\bibitem{6}
Farhadian, R. \emph{A probabilistic proof of a recursion formula for sums of integer powers,} Integers \textbf{23} (2023), Paper No. A88. 
\bibitem{7}
Hu, X.; Zhong, Y. \emph{A probabilistic proof of a recursion formula for sums of powers,} Amer. Math. Monthly \textbf{127} (2020), no. 2, 166-168.
\bibitem{8}
Kim, D. S.; Kim, T. \emph{A note on a new type of degenerate Bernoulli numbers,} Russ. J. Math. Phys. \textbf{27} (2020), no. 2, 227-235.
\bibitem{9}
Kim, G.; Kim, B.; Choi, J. \emph{The $DC$-algorithm for computing sums of powers of consecutive integers and Bernoulli numbers,} Adv. Stud. Contemp. Math. (Kyungshang) \textbf{17} (2008), no. 2, 137-145.
\bibitem{10}
Kim, T. \emph{Sums of powers of consecutive $q$-integers,} Adv. Stud. Contemp. Math. (Kyungshang) \textbf{9} (2004), no. 1, 15-18.
\bibitem{11}
Kim, T.; Kim, D. S. \emph{Generalization of Spivey's recurrence relation,} Russ. J. Math. Phys. \textbf{31} (2024), no. 2, 218-226.
\bibitem{12}
Kim, T.; Kim, D. S. \emph{Probabilistic degenerate Bell polynomials associated with random variables,} Russ. J. Math. Phys. \textbf{30} (2023), no. 4, 528-542.
\bibitem{13}
Kim, T.; Kim, D. S. \emph{A new approach to fully degenerate Bernoulli numbers and polynomials,} Filomat \textbf{37} (2023), no. 7, 2269-2278.
\bibitem{14}
Kim, T.; Kim, D. S.; Kim, H. K. \emph{Some identities related to degenerate
Stirling numbers of the second kind,} Demonstr. Math. \textbf{55} (2022), 812–821.
\bibitem{15}
Kim, T.; Rim, S.-H.; Simsek, Y. \emph{A note on the alternating sums of powers of consecutive $q$-integers,} Adv. Stud. Contemp. Math. (Kyungshang) \textbf{13} (2006), no. 2, 159-164.
\bibitem{16}
Luo, L.; Kim, T.; Kim, D. S.; Ma, Y. \emph{Probabilistic degenerate Bernoulli and degenerate Euler polynomials,} Math. Comput. Model. Dyn. Syst. \textbf{30} (2024), no. 1, 342-363.
\bibitem{17}
Park, J.-W.; Pyo, S.-S. \emph{A note on degenerate Bernoulli polynomials arising from umbral calculus,} Adv. Stud. Contemp. Math. (Kyungshang) \textbf{32} (2022), no. 4, 509-525.
\bibitem{18}
Simsek, Y. \emph{Identities and relations related to combinatorial numbers and polynomials,} Proc. Jangjeon Math. Soc. \textbf{20} (2017), no. 1, 127-135.
\bibitem{19}
Spivey, M. Z. \emph{A combinatorial view of sums of powers,} Math. Mag. \textbf{94} (2021), no. 2, 125-131. 
\bibitem{20}
Xu, R.; Kim, T.; Kim, D. S.; Ma, Y. \emph{Probabilistic degenerate Fubini polynomials associated with random variables,} J. Nonlinear Math. Phys. \textbf{31} (2024), no. 1, Paper No. 47. 
\end{thebibliography}
\end{document}